\input amstex
\documentstyle{amsppt}
%----------------------------------------------------------------
% Title:     Transfinite normal and composition series of modules.
% Author:    Ruslan Sharipov
% Comments:  AmSTeX, 7 pages, amsppt style
% MSC-class: 16D70
%----------------------------------------------------------------
%           Replacement for output macro definition
%
\catcode`@=11
\redefine\output@{%
  \def\break{\penalty-\@M}\let\par\endgraf
  \ifodd\pageno\global\hoffset=105pt\else\global\hoffset=8pt\fi  
  \shipout\vbox{%
    \ifplain@
      \let\makeheadline\relax \let\makefootline\relax
    \else
      \iffirstpage@ \global\firstpage@false
        \let\rightheadline\frheadline
        \let\leftheadline\flheadline
      \else
        \ifrunheads@ %\let\makefootline\relax
        \else \let\makeheadline\relax
        \fi
      \fi
    \fi
    \makeheadline \pagebody \makefootline}%
  \advancepageno \ifnum\outputpenalty>-\@MM\else\dosupereject\fi
}
\catcode`\@=\active
%----------------------------------------------------------------
\nopagenumbers
\def\negskp{\hskip -2pt}

\def\End{\operatorname{End}}
\def\Ker{\operatorname{Ker}}
\def\compos{\,\raise 1pt\hbox{$\sssize\circ$} \,}
\def\Cl{\operatorname{Cl}}

\def\blue#1{#1}

\catcode`#=11\def\diez{#}\catcode`#=6
\catcode`_=11\def\podcherkivanie{_}\catcode`_=8
%\catcode`~=11\def\volna{~}\catcode`~=\active
\def\mycite#1{\cite{\blue{#1}}\immediate\special{ps:
     ShrHPSdict begin /ShrBORDERthickness 0 def}}

\def\mytag#1{%
    \tag#1}
\def\mythetag#1{\thetag{\blue{#1}}\immediate\special{ps:
     ShrHPSdict begin /ShrBORDERthickness 0 def}}
\def\myrefno#1{\no#1}
\def\myhref#1#2{\blue{#2}\immediate\special{ps:
     ShrHPSdict begin /ShrBORDERthickness 0 def}}
\def\myEarXivlink{\myhref{http://arXiv.org}{http:/\negskp/arXiv.org}}

\def\mytheorem#1{\csname proclaim\endcsname{Theorem #1}}
\def\mytheoremwithtitle#1#2{\csname proclaim\endcsname{Theorem #1#2}}
\def\mythetheorem#1{\blue{#1}\immediate\special{ps:
     ShrHPSdict begin /ShrBORDERthickness 0 def}}
\def\mylemma#1{\csname proclaim\endcsname{Lemma #1}}
\def\mylemmawithtitle#1#2{\csname proclaim\endcsname{Lemma #1#2}}
\def\mythelemma#1{\blue{#1}\immediate\special{ps:
     ShrHPSdict begin /ShrBORDERthickness 0 def}}
\def\mycorollary#1{\csname proclaim\endcsname{Corollary #1}}

\def\mydefinition#1{\definition{Definition #1}}
\def\mythedefinition#1{\blue{#1}\immediate\special{ps:
     ShrHPSdict begin /ShrBORDERthickness 0 def}}

%----------------------------------------------------------------
% Cyrillic fonts definition
%\font\eightcyr=wncyr8
%----------------------------------------------------------------
\pagewidth{360pt}
\pageheight{606pt}
\topmatter
\title
Transfinite normal and composition series of modules.
\endtitle
\author
R.~A.~Sharipov
\endauthor
\address 5 Rabochaya street, 450003 Ufa, Russia\newline
\vphantom{a}\kern 12pt Cell Phone: +7(917)476 93 48
\endaddress
\email \vtop to 30pt{\hsize=280pt\noindent
\myhref{mailto:r-sharipov\@mail.ru}
{r-sharipov\@mail.ru}\newline
\myhref{mailto:R\podcherkivanie Sharipov\@ic.bashedu.ru}
{R\_\hskip 1pt Sharipov\@ic.bashedu.ru}\vss}
\endemail
\urladdr
\vtop to 20pt{\hsize=280pt\noindent
\myhref{http://www.freetextbooks.narod.ru}
{http:/\negskp/www.freetextbooks.narod.ru}\newline
\myhref{http://sovlit2.narod.ru/}
{http:/\negskp/sovlit2.narod.ru}\vss}
\endurladdr
\abstract
    Normal and composition series of modules enumerated by ordinal numbers 
are studied. The Jordan-H\"older theorem for them is discussed. 
\endabstract
\endtopmatter
%\loadbold
%\loadeufb
\TagsOnRight
\document

\head
1. Introduction and some preliminaries.
\endhead
     Transfinite normal and composition series of groups were studied
in \mycite{1}. They generalize classical normal and composition series 
which are finitely long. In this paper we reproduce the results of 
\mycite{1} in the case of left modules over associative rings and 
associative algebras. 
\mydefinition{1.1} Let $A$ be an associative ring and let $V$ be an
additive Abelian group. The Abelian group $V$ is called a left module 
over the ring $A$ (or a left $A$-module) if some homomorphism of rings 
$$
\varphi\!:\,A\to\End(V) 
$$
is given and fixed. Here $\End(V)$ is the ring of endomorphisms of the 
additive Abelian group $V$. 
\enddefinition
\mydefinition{1.2} Let $A$ be an associative algebra over some field $\Bbb K$
and let $V$ be a linear vector space over the same field $\Bbb K$. 
The linear vector space $V$ is called a left module over the algebra $A$ 
(left $A$-module) if some homomorphism of algebras
$$
\varphi\!:\,A\to\End(V) 
$$
is given and fixed. Here $\End(V)$ is the algebra of endomorphisms of the 
linear vector space $V$. 
\enddefinition
     In both cases each element $a\in A$ produces an operator $\varphi(a)$
that can act upon any element $v\in V$. The result of applying $\varphi(a)$
to $v$ is denoted as
$$
\hskip -2em
u=\varphi(a)(v).
\mytag{1.1}
$$
In some cases the formula \mythetag{1.1} is written as follows:
$$
\hskip -2em
u=\varphi(a)v.
\mytag{1.2}
$$
The endomorphism $\varphi(a)\in\End(V)$ in this formula is presented as a
left multiplier for the element $v\in V$. In some cases the symbol $\varphi$ 
is also omitted:
$$
\hskip -2em
u=a\,v.
\mytag{1.3}
$$
The formula \mythetag{1.3} is a purely algebraic version of the formulas
\mythetag{1.1} and \mythetag{1.2}. The action of $\varphi(a)$ upon $v$ here
is written as a multiplication of $v$ on the left by the element $a\in A$. 
This is the reason for which $V$ is called a left $A$-module. Using 
\mythetag{1.3}, the definitions~\mythedefinition{1.1} and \mythedefinition{1.2} 
can reformulated as follows. 
\mydefinition{1.3} A left module over an associative ring $A$ is an additive 
Abelian group $V$ equipped with an auxiliary operation of left multiplication 
by elements of the ring $A$ such that the following identities are fulfilled: 
\roster
\rosteritemwd=0pt
\item"1)" \ $a\,(v_1+v_2)=a\,v_1+a\,v_2$ for all $a\in A$ and for all $v_1,v_2\in V$;
\item"2)" \ $(a_1+a_2)\,v=a_1\,v+a_2\,v$ for all $a_1,a_2\in A$ and for all $v\in V$;
\item"3)" \ $(a_1\,a_2)\,v=a_1\,(a_2\,v)$ for all $a_1,a_2\in A$ and for all $v\in V$.
\endroster
\enddefinition
\mydefinition{1.4} A left module over an associative $\Bbb K$-algebra $A$ is a
linear vector space  $V$ over the field $\Bbb K$ equipped with an auxiliary operation 
of left multiplication by elements of the algebra $A$ such that the following 
identities are fulfilled: 
\roster
\rosteritemwd=0pt
\item"1)" \ $a\,(v_1+v_2)=a\,v_1+a\,v_2$ for all $a\in A$ and for all $v_1,v_2\in V$;
\item"2)" \ $(a_1+a_2)\,v=a_1\,v+a_2\,v$ for all $a_1,a_2\in A$ and for all $v\in V$;
\item"3)" \ $(a_1\,a_2)\,v=a_1\,(a_2\,v)$ for all $a_1,a_2\in A$ and for all $v\in V$;
\item"4)" \ $(k\,a)\,v=a\,(k\,v)=k\,(a\,v)$ for all $k\in\Bbb K$, $a\in A$, and $v\in V$.
\endroster
\enddefinition
     Right $A$-modules are similar to left ones. Here elements $v\in V$
are multiplied by elements $a\in A$ on the right. However, each right
$A$-module can be treated as a left $A^\bullet$-module, where $A^\bullet$ 
is the opposite ring (opposite algebra) for $A$ (see Chapter~5 in \mycite{2}).
For this reason below we consider left modules only.
\head
2. Submodules and factormodules. 
\endhead
\mydefinition{2.1} Let $V$ be a left $A$-module, where $A$ is a ring. A subset 
$W\subseteq V$ is called a submodule of $V$ if it is closed with respect to
the inversion, with respect to the addition, and with respect to the multiplication 
by elements of the ring $A$:
\roster
\rosteritemwd=0pt
\item"1)" \ $v\in W$ implies $-v\in W$;
\item"2)" \ $v_1\in W$ and $v_2\in W$ imply $(v_1+v_2)\in W$;
\item"3)" \ $v\in W$ and $a\in A$ imply $(a\,v)\in W$.
\endroster
In other words, a submodule $W$ is a subgroup of the additive group $V$ invariant
with respect to the endomorphisms $\varphi(a)$ for all $a\in A$. 
\enddefinition
\mydefinition{2.2} Let $V$ be a left $A$-module, where $A$ is some $\Bbb K$-algebra. 
A subset $W\subseteq V$ is called a submodule of $V$ if it is closed with respect 
to the addition and with respect to the multiplication by elements of the field 
$\Bbb K$ and the algebra $A$:
\roster
\rosteritemwd=0pt
\item"1)" \ $v\in W$ and $k\in\Bbb K$ imply $(k\,v)\in W$;
\item"2)" \ $v_1\in W$ and $v_2\in W$ imply $(v_1+v_2)\in W$;
\item"3)" \ $v\in W$ and $a\in A$ imply $(a\,v)\in W$.
\endroster
In other words, a submodule $W$ is a subspace of the vector space $V$ invariant
with respect to the endomorphisms $\varphi(a)$ for all $a\in A$. 
\enddefinition
     In both cases a submodule $W\subseteq V$ inherits the structure of a
left $A$-module from $V$. Each submodule $W$ is associated with the corresponding 
factorset $V/W$. The factorset $V/W$ is composed by cosets
$$
\hskip -2em
\Cl_W(v)=\{u\in V\!:\,u=v+w\text{\ \ for some \ }w\in W\}.
\mytag{2.1}
$$
The element $v$ in \mythetag{2.1} is a representative of the coset $\Cl_W(v)$.
Any element $v\in\Cl_W(v)$ can be chosen as its representative. In both cases 
(where $A$ is a ring or where $A$ is an algebra) the factorset $V/W$ inherits 
the structure of a left $A$-module from $V$. For this reason it is called a 
factormodule. Algebraic operations with cosets are given by the formulas
\roster
\rosteritemwd=0pt
\item"1)" \ $\Cl_W(v_1)+\Cl_W(v_2)=\Cl_W(v_1+v_2)$ for any $v_1,v_2\in V$;
\item"2)" \ $a\Cl_W(v)=\Cl_W(a\,v)$ for any $a\in A$ and for any $v\in V$.
\endroster
\noindent
If $A$ is a $\Bbb K$-algebra, then
\roster
\rosteritemwd=0pt
\item"3)" \ $k\Cl_W(v)=\Cl_W(k\,v)$ for any $k\in\Bbb K$ and for any $v\in V$.
\endroster
\head
3. Transfinite normal and composition series.
\endhead
\mydefinition{3.1} Let $V$ be a left $A$-module. A transfinite sequence of its
submodules 
$$
\hskip -2em
\{0\}=V_1\varsubsetneq V_2\varsubsetneq\ldots
\varsubsetneq V_n=V
\mytag{3.1}
$$
is called a transfinite normal series for $V$ if  
$V_\alpha=\dsize\bigcup_{\beta<\alpha}V_\beta$ for any limit ordinal 
$\alpha\leqslant n$.
\enddefinition
     The term ``normal series'' in the above definition comes from the group 
theory where each normal series $\{1\}=G_1\varsubsetneq G_2\varsubsetneq\ldots
\varsubsetneq G_n=G$ of a group $G$ should be composed by its subgroups such that 
$G_i$ is a normal subgroup of $G_{i+1}$. Otherwise we would not be able to 
build the factorgroup $G_{i+1}/G_i$. In the case of modules each submodule 
$V_i$ with $i<n$ produces the factormodule $V_{i+1}/V_i$. So the term ``normal 
series'' here is vestigial, it is used for ear comfort only. 
\mydefinition{3.2} A module $V$ is called hypertranssimple if it has 
no normal series (neither finite nor transfinite) other than trivial one 
$\{0\}=V_1\varsubsetneq V_2=V$.
\enddefinition
\mydefinition{3.3} A transfinite normal series \mythetag{3.1} of a module $V$ 
is called a transfinite composition series of $V$ if for each ordinal number
$i<n$ the factormodule $V_{i+1}/V_i$ is hypertranssimple.
\enddefinition
     The concept of hypertranssimplicity is nontrivial in the case of groups.
In the case of modules each nontrivial submodule $W$ of $V$ produces the nontrivial
normal series $\{0\}=V_1\varsubsetneq V_2\varsubsetneq V_3=V$, where $V_2=W$. For
this reason the concept of hypertranssimplicity of modules reduces to the standard 
concept of simplicity. The definitions~\mythedefinition{3.2} and 
\mythedefinition{3.3} then are reformulated as follows.
\mydefinition{3.4} A module $V$ is called simple if it has no submodules other
than trivial ones $V_1=\{0\}$ and $V_2=V$.
\enddefinition
\mydefinition{3.5} A transfinite normal series \mythetag{3.1} of a module $V$ 
is called a transfinite composition series of $V$ if for each ordinal number 
$i<n$ the corresponding factormodule $V_{i+1}/V_i$ is simple.
\enddefinition
\head
4. Intersections and sums of submodules.
\endhead
     The intersection of two or more submodules of a given module $V$ is again 
a submodule of $V$. As for unions, the union of submodules in general case is 
not a submodule. Unions of submodules are used to define their sums. 
\mydefinition{4.1} Let $U_i$ with $i\in I$ be submodules of some module $V$.
The submodule $U$ generated by the union of submodules $U_i$ is called their sum:
$$
\hskip -2em
U=\sum_{i\in I}U_i=\Big<\bigcup_{i\in I}U_i\Big>.
\mytag{4.1}
$$
\enddefinition
\noindent If the number of submodules in \mythetag{4.1} is finite, one can use
the following notations:
$$
\hskip -2em
U=U_1+\ldots+U_n=\langle\kern 1pt U_1\cup\ldots\cup\kern 1pt U_n\rangle.
\mytag{4.2}
$$
In both cases \mythetag{4.1} and \mythetag{4.2} each element $u\in U$ is presented
as a finite sum
$$
u=u_{i_1}+\ldots+u_{i_s}\text{, \ where \ }u_{i_r}\in U_{i_r}
\text{\ \ and \ }i_r\in I\text{\ \ for all \ }r=1,\,\ldots,\,s.
\quad
\mytag{4.3}
$$
\mydefinition{4.2} The sum of submodules \mythetag{4.1} is called a direct sum 
if for each element $u\in U$ its presentation \mythetag{4.3} is unique. 
\enddefinition
\mylemmawithtitle{4.1}{ (Zassenhaus)} Let\/ $\tilde U$ and $\tilde W$ be 
submodules of some left $A$-module and let $U$ and $W$ be submodules
of\/ $\tilde U$ and $\tilde W$ respectively. Then 
$$
\hskip -2em
(U+(\tilde U\cap\tilde W))/(U+(\tilde U\cap W))\cong 
(W+(\tilde W\cap\tilde U))/(W+(\tilde W\cap U)).
\mytag{4.4}
$$
\endproclaim
     The lemma~\mythelemma{4.1} is also known as the butterfly lemma. 
Typically the butterfly lemma is formulated for groups (see \S\,3 of Chapter~\uppercase\expandafter{\romannumeral 1} in \mycite{3}).
However, its proof can be easily adapted for the case of modules.
\demo{Proof} Let's denote $M=U+(\tilde U\cap W)$ and $N=W+(\tilde W
\cap U)$. Elements of the factormodule in the left hand side of the formula
\mythetag{4.4} are cosets of the form
$$
\hskip -2em
\Cl_M(u+a)\text{, \ where \ }u\in U\text{\ \ and \ }a\in\tilde U\cap\tilde W.
$$
Note that $U\subseteq M$. Therefore $\Cl_M(u+a)=\Cl_M(a)$ which means that
each element of the factormodule $(U+(\tilde U\cap\tilde W))/M$ is represented
by some element $a\in\tilde U\cap\tilde W$. Thus we have a surjective 
homomorphism of modules
$$
\hskip -2em
\varphi\!:\,\tilde U\cap\tilde W\longrightarrow(U+(\tilde U\cap\tilde W))/M.
\mytag{4.5}
$$
Repeating the above arguments for the factormodule in the right hand side of 
the formula \mythetag{4.4}, we get another surjective homomorphism of modules
$$
\hskip -2em
\psi\!:\,\tilde U\cap\tilde W\longrightarrow(W+(\tilde W\cap\tilde U))/N.
\mytag{4.6}
$$\par
     The rest is to prove that $\Ker\varphi=\Ker\psi$. Assume that $a\in 
\Ker\varphi$. In this case $a\in\tilde U\cap\tilde W$ and $a\in M$. The 
inclusion $a\in M$ means that $a=u+w$, where $u\in U$ and $w\in\tilde U
\cap W$. \pagebreak Since $\tilde U\cap W\subseteq\tilde U\cap\tilde W$, 
from $u=a-w$ we derive $u\in\tilde U\cap\tilde W$. On the other hand $u\in U$. 
Hence $u\in U\cap(\tilde U\cap\tilde W)$, which means $u\in U\cap\tilde W$. 
Thus we have proved that each element $a\in\Ker\varphi$ is presented as a sum
$$
\hskip -2em
a=u+w\text{, \ where \ }u\in\tilde W\cap U\text{\ \ and \ }w\in\tilde U\cap W.
\mytag{4.7}
$$
Conversely, it is easy to see that the presentation \mythetag{4.7} leads to 
$a\in\tilde W\cap\tilde U$ and $a\in M$, i\.\,e\. $a\in\Ker\varphi$. For this 
reason we have $\Ker\varphi=(\tilde W\cap U)+(\tilde U\cap W)$. The equality
$\Ker\psi=(\tilde W\cap U)+(\tilde U\cap W)$ is proved similarly. Now the
formula \mythetag{4.4} is immediate from $\Ker\varphi=\Ker\psi$ due to the
surjectivity of the homomorphisms \mythetag{4.5} and \mythetag{4.6}. The
butterfly lemma~\mythelemma{4.1} is proved.
\qed\enddemo 
\head
5. The Jordan-H\"older theorem.
\endhead
\mydefinition{5.1} A transfinite normal series $\{0\}=\tilde V_1\varsubsetneq 
\tilde V_2\varsubsetneq\ldots\varsubsetneq\tilde V_p=G$ is called a refinement
for a transfinite normal series $\{0\}=V_1\varsubsetneq V_2\varsubsetneq\ldots
\varsubsetneq V_n=G$ if each submodule $V_i$ coincides with some submodule
$\tilde V_j$.
\enddefinition
\mydefinition{5.2} Two transfinite normal series $\{0\}=V_1\varsubsetneq V_2
\varsubsetneq\ldots\varsubsetneq V_n=V$ and $\{0\}=W_1\varsubsetneq W_2
\varsubsetneq\ldots\varsubsetneq W_m=V$ of a module $V$ are called isomorphic 
if there is a one-to-one mapping that associates each ordinal number $i<n$
with some ordinal number $j<m$ in such a way that 
$V_{i+1}/V_i\cong W_{j+1}/W_j$.
\enddefinition
\mytheorem{5.1} Arbitrary two transfinite normal series of a left $A$-module 
$V$ have isomorphic refinements. 
\endproclaim
\mylemma{5.1} If\/ $\{1\}=G_1\varsubsetneq\ldots \varsubsetneq G_n=G$ 
is a transfinite composition series of a group $G$, then it has no refinements 
different from itself. 
\endproclaim
\mytheoremwithtitle{5.2}{ (Jordan-H\"older)}Any two transfinite composition 
series of a left $A$-module $V$ are isomorphic. 
\endproclaim
     The theorem~\mythetheorem{5.2} is immediate from the 
theorem~\mythetheorem{5.1} and the lemma~\mythelemma{5.1}. As for the
theorem~\mythetheorem{5.1} and the lemma~\mythelemma{5.1}, their proof is
quite similar to the proof of the theorem~{3.1} and the lemma~{3.10} in
\mycite{1}. The algebraic part of this proof is based on the Zassenhaus
butterfly lemma. Its version for modules is given above (see 
lemma~\mythelemma{4.1}). The other part of the proof deals with indexing 
sets and ordinal numbers, not with algebraic structures. For this reason it 
does not differ in the case of groups and in the case of modules.
\head
6. External direct sums. 
\endhead
     Sums and direct sums introduced in the definitions~\mythedefinition{4.1}
and \mythedefinition{4.2} are internal ones. They are formed by submodules
of a given module. External direct sums are formed by separate modules which
are not necessarily enclosed in a given module. 
\mydefinition{6.1} Let $V_i$ be left $A$-modules enumerated by elements $i\in I$
of some indexing set $I$. Finite formal sums of the form
$$
v=v_{i_1}+\ldots+v_{i_s}\text{, \ where \ }v_{i_r}\in V_{i_r}
\text{\ \ and \ }i_r\in I\text{\ \ for all \ }r=1,\,\ldots,\,s,
\quad
\mytag{6.1}
$$
constitute a left $A$-module $V$ which is called the direct sum of the 
modules $V_i$. 
\enddefinition
     Once the external direct sum $V$ is constructed, we find that it
comprises submodules $U_i$ isomorphic to the initial modules $V_i$. 
Indeed, we can set $s=1$ in \mythetag{6.1}. Formal sums \mythetag{6.1} 
with exactly one summand $v=v_i$, where $v_i\in V_i$, constitute a submodule
$U_i$ of $V$ isomorphic to the module $V_i$. For this reason the constructions 
of internal and external direct sums are the same in essential.\par
      According to the well-known Zermelo theorem (see Appendix~2 in
\mycite{4}), every set $I$ can be well ordered and then associated with some
ordinal number $n$ (see Proposition~3.8 in Appendix~3 of \mycite{4}). Therefore
the external direct sum $V$ in the above definition~\mythedefinition{6.1} can 
be written as
$$
\hskip -2em
V=\bigoplus_{\alpha<n}V_\alpha.
\mytag{6.2}
$$
Relying on \mythetag{6.2}, for each $i\leqslant n$ we introduce the following 
submodule of $V$:
$$
\hskip -2em
W_i=\bigoplus_{\alpha<i}V_\alpha.
\mytag{6.3}
$$
The submodules \mythetag{6.3} constitute a transfinite normal series of submodules
for $V$:
$$
\hskip -2em
\{0\}=W_1\varsubsetneq W_2\varsubsetneq\ldots\varsubsetneq W_n=V.
\mytag{6.4}
$$
The factormodules of the sequence \mythetag{6.4} are isomorphic to the modules
$V_i$, namely we have $W_{i+1}/W_i\cong V_i$. If the modules $V_i$ are simple,
then \mythetag{6.4} is a transfinite composition series (see 
Definitions~\mythedefinition{3.4} and \mythedefinition{3.5}).\par
     {\bf A remark}. The formulas \mythetag{6.2}, \mythetag{6.3}, and \mythetag{6.4}
are equally applicable for internal and external direct sums. For this reason, 
applying the Jordan-H\"older theorem~\mythetheorem{5.2}, we get the following result. 
\mytheorem{6.1} If a module $V$ is presented as a direct sum of its simple submodules
$V_i$, then these submodules are unique up to the isomorphism and some permutation of
their order in the direct sum. 
\endproclaim
     There is a special case of the external direct sum \mythetag{6.2}. Assume that
$U$ is some simple left $A$-module. Let's replicate this module into multiple copies
and denote these copies through $V_\alpha$. Then $V_\alpha\cong U$. In this case the
module $V$ in \mythetag{6.2} is denoted through $N\,U$, where $N=|n|$ is the cardinality 
of the ordinal number $n$. Such a notation is motivated by the following theorem.
\mytheorem{6.2} Let $U$ be a simple left $A$-module and let $N\,U$ and $M\,U$ be
two external direct sums of the form \mythetag{6.2} built by the copies of the
module $U$:
$$
\xalignat 2
&N\,U=\bigoplus_{\alpha<n}U,
&&M\,U=\bigoplus_{\alpha<m}U.
\endxalignat
$$
Then $N\,U$ is isomorphic to $M\,U$ if and only if\/ $N=M$, i\.\,e\. if\/ $|n|=|m|$. 
\endproclaim
The theorem~\mythetheorem{6.2} is easily derived from the theorem~\mythetheorem{6.1}.
\head
7. Concluding remarks.
\endhead
    The results of this paper are rather obvious and are known to the algebraists
community. However, they are dispersed in various books as preliminaries to more
special theories. Treated as obvious, these results are usually not equipped with
explicit proofs and even with explicit statements. We gather them in this paper 
for referential purposes. 

\enddocument
\end